\documentclass[conference]{IEEEtran}
\usepackage{amsmath,amssymb,amsthm,epsfig,dsfont,color,subfigure,empheq,enumerate,url}
\usepackage{algorithmic,algorithm,balance}
\usepackage{textcomp}

\newcommand \bone{\mathbf{1}}

\newcommand \bc{\mathbf{c}}

\newcommand \bdell{\boldsymbol{\ell}} 

\newcommand \bn{\mathbf{n}}

\newcommand \bt{\mathbf{t}}

\newcommand \bx{\mathbf{x}}
\newcommand \by{\mathbf{y}}

\newcommand \bpi{\boldsymbol{\pi}}

\newcommand \mcC{\mathcal{C}}

\newcommand \mcE{\mathcal{E}}

\newcommand \mcG{\mathcal{G}}

\newcommand \mcL{\mathcal{L}}

\newcommand \mcN{\mathcal{N}}

\newcommand \mcP{\mathcal{P}}

\newcommand \mcS{\mathcal{S}}

\allowdisplaybreaks

\begin{document}

\title{Optimal Distribution System Restoration with \\ Microgrids and Distributed Generators}
\author{\IEEEauthorblockN{Manish Kumar Singh, Vassilis Kekatos, and Chen-Ching Liu}
	\IEEEauthorblockA{Emails:\{manishks,kekatos,ccliu\}@vt.edu}
Bradley Dept. of Electrical \& Computer Engnr., Virginia Tech, Blacksburg, VA 24061, USA}

\maketitle

\begin{abstract}
Increasing emphasis on reliability and resiliency call for advanced distribution system restoration (DSR). The integration of grid sensors, remote controls, and distributed generators (DG) brings about exciting opportunities in DSR. In this context, this work considers the task of single-step restoration of a single-phase power distribution system. Different from existing works, the devised restoration scheme achieves optimal formation of islands without heuristically pre-identifying reference buses. It further facilitates multiple DGs running within the same island, and establishes a coordination hierarchy in terms of their PV/PQ operation modes. Generators without black-start capability are guaranteed to remain connected to a black-start DG or a substation. The proposed scheme models remotely-controlled voltage regulators exactly, and integrates them in the restoration process. Numerical tests on a modified IEEE 37-bus feeder demonstrate that the proposed mixed-integer linear program (MILP) takes less than four seconds to handle random outages of 1--5 lines. The scalability of this novel MILP formulation can be attributed to the unique use of cycles and paths on the grid infrastructure graph; the McCormick linearization technique; and an approximate power flow model. 
\end{abstract}


\section{Introduction}
Outages in power distribution systems are inevitable and could range from single-line outages caused by faults to widespread outages due to extreme events~\cite{jianhui2017proc}, \cite{ccliu2018coordinating}. Extreme events in the form of natural disasters, accidents, or cyber attacks, could result in a tremendous loss of distribution infrastructure \cite{baldick2016resilience}, \cite{jianhui2018sequential}. Meticulously designing the response to outages could significantly improve system resiliency~\cite{ton2015magazine}. Conventionally, distribution system restoration was predominantly manual, and was based on trouble calls, the operator's prior experience, and field search crews. Currently, the rampant deployment of smart meters, grid sensors, and controlled switches, offers improved situational awareness and remote control~\cite{abbey2014magazine}. To expedite grid restoration, many efforts have been put towards its automation based on available resources. 

The task of distribution system restoration (DSR) is initiated after the post-outage status of the distribution grid has been assessed. Operators traditionally resort to network reconfiguration schemes to limit the impact of an outage while satisfying operational and design constraints. The DSR task is typically formulated as a combinatorially complex non-linear minimization. Traditionally, it has been approached using dynamic programming \cite{heydt2008DP}; expert systems \cite{ccliu1988expert}; fuzzy logic \cite{slee2006fuzzy}; genetic algorithms \cite{genetic2004restoration}; and mixed-integer non-linear programming \cite{jignesh2007minlp}. Nonetheless, with the advent of distributed generators (DG) and microgrids, new challenges and opportunities have been introduced in the DSR problem. 

Distributed generators and microgrids could enable islanded operation, thus improving resiliency against extreme events. The coordinated operation of heterogeneous DGs introduces different operational and control requirements~\cite{iravani2006MG}. Although several recent works deal with DGs and microgrids~\cite{jianhui2018sequential}, \cite{ccliu2014spanning}, \cite{jianhui2015MINLP}, \cite{franco2013milp}; they all presume that each DG features black-start capability and/or preclude running multiple DGs on the same island. The former does not hold for solar DGs without energy storage. The latter over-simplifies the operational capabilities of DGs and thus constitutes a restriction of the actual DSR task. Albeit \cite{ccliu2018coordinating} allows for multiple DGs operating on the same island, their control mode is decided through a suboptimal two-stage process. In a nutshell, a realistic coordination of DGs and microgrids for DSR remains largely under-explored.



This work puts forth a novel system restoration scheme with three major improvements over existing alternatives: \emph{i)} Our DSR scheme finds the optimal formation of islands in a single stage, different from prior works that first identify reference generators and then build islands around them; \emph{ii)} It further allows for multiple (non)-black-start DGs running on the same island and decides their optimal coordination; and \emph{iii)} It devises an exact model for voltage regulators. Through the novel use of cycles and paths on the grid infrastructure graph, and by leveraging the McCormick linearization and an approximate grid model, our optimal DSR task can be posed as an MILP, which scales well on a moderately-sized feeder.

\section{ Preliminaries}\label{sec:pre}
Before defining the DSR task, we review some preliminaries from graph theory and the McCormick linearization. Consider an undirected graph $\mcG := (\mcN,\mcE)$ with $\mcN$ being its node set and $\mcE$ its edge set. The graph is connected if there exists a sequence of adjacent edges between any two of its nodes. A \emph{path} from node $i$ to $j$ is defined as the sequence of edges $\mcP_{ij}$ starting at node $i$ and terminating at node $j$. 
 A \emph{cycle} is a sequence of adjacent edges without repetition that starts and ends at the same node. A \emph{tree} is a connected graph with no cycles. If every edge $e\in\mcE$ is assigned a direction, the obtained graph is termed \emph{directed}. 

As a brief review, the McCormick linearization is a widely used technique for handling products of optimization variables $x_1x_2\cdot x_N$ by their linear convex envelopes~\cite{McCormick1976}. Since this relaxation is not necessarily exact, there is a rich literature on tightening approaches; see for example \cite{harsha2016mccormick} and references therein. In fact, the McCormick linearization becomes exact for the special case of bilinear terms involving at least one binary variable: Consider the constraint $z=xy$, according to which the variable $z$ equals the the binary variable $x\in\{0,1\}$ times the continuous variable $y$. If $y$ is constrained within $y\in[\underline{y},\bar{y}]$, the constraint $z=xy$ can be equivalently expressed by the four linear inequality constraints
\begin{subequations}\label{eq:MC}
\begin{align}
x\underline{y}&\leq z\leq x\bar{y},\label{seq:MC1}\\
y+(x-1)\bar{y}&\leq z\leq y +(x-1)\underline{y}\label{seq:MC2}.
\end{align}
\end{subequations}
The equivalence can be readily verified by observing that for $x=1$, constraint \eqref{seq:MC2} yields $z=y$ and \eqref{seq:MC1} holds trivially. When $x=0$, both \eqref{seq:MC1} and \eqref{seq:MC2} yield $z=0$. Combining the two cases provides $z=xy$ indeed. Henceforth, all bilinear products of binary and (bounded) continuous variables could be handled by the McCormick linearization of \eqref{eq:MC}.

\section{Problem Formulation}\label{sec:problem}
During an outage, protection devices isolate certain parts of a distribution system including the faulty elements. While replacing faulty elements may be time-consuming, remotely-controlled switches could reconfigure the system to alleviate the outage effect. Given the post-fault status, the \emph{single-step DSR} task finds the grid topology that minimizes the outage impact while complying with operational constraints. The presumption is that the system transitions instantaneously from the post-fault to the final condition. In practice, this transition is implemented via a sequence of reconfiguration steps involving one control action at a time. Due to space limitations, here we consider the single-step DSR task on a single-phase grid model. 

A distribution system can be represented by a graph $\mcG:=(\mcN,\mcE)$. Its nodes are indexed by $i\in\mcN:=\{0,\dots,N\}$ correspond to buses, and its edges $\mcE$ to distribution lines, switches, and voltage regulators. An edge running between nodes $i$ and $j$ is assigned an arbitrary direction, and is denoted as $e:(i,j)$ or $e:(j,i)\in\mcE$. Although multiple islands may be formed by opening switches, graph $\mcG$ is connected since the distribution system is structurally connected.

\subsection{Nodal Variables and Constraints}\label{subsec:nodes}
Each bus $i\in\mcN$ hosts at most one generator or load. This is without loss of generality since a bus with multiple loads can be modeled as a set of single-load buses, all connected by non-switchable zero-impedance lines. Moreover, to ensure that all substations remain disconnected from each other, they are combined into a single root node indexed by $0$ as in~\cite{ccliu2014spanning}. The power limit of a substation can be imposed as a limit on the line connecting the substation with its feeder. To simplify the exposition, on-load tap changers (OLTCs) are ignored and all substations are assumed to operate at the nominal voltage.

\begin{figure}[t]
	\centering
	\includegraphics[scale=0.3]{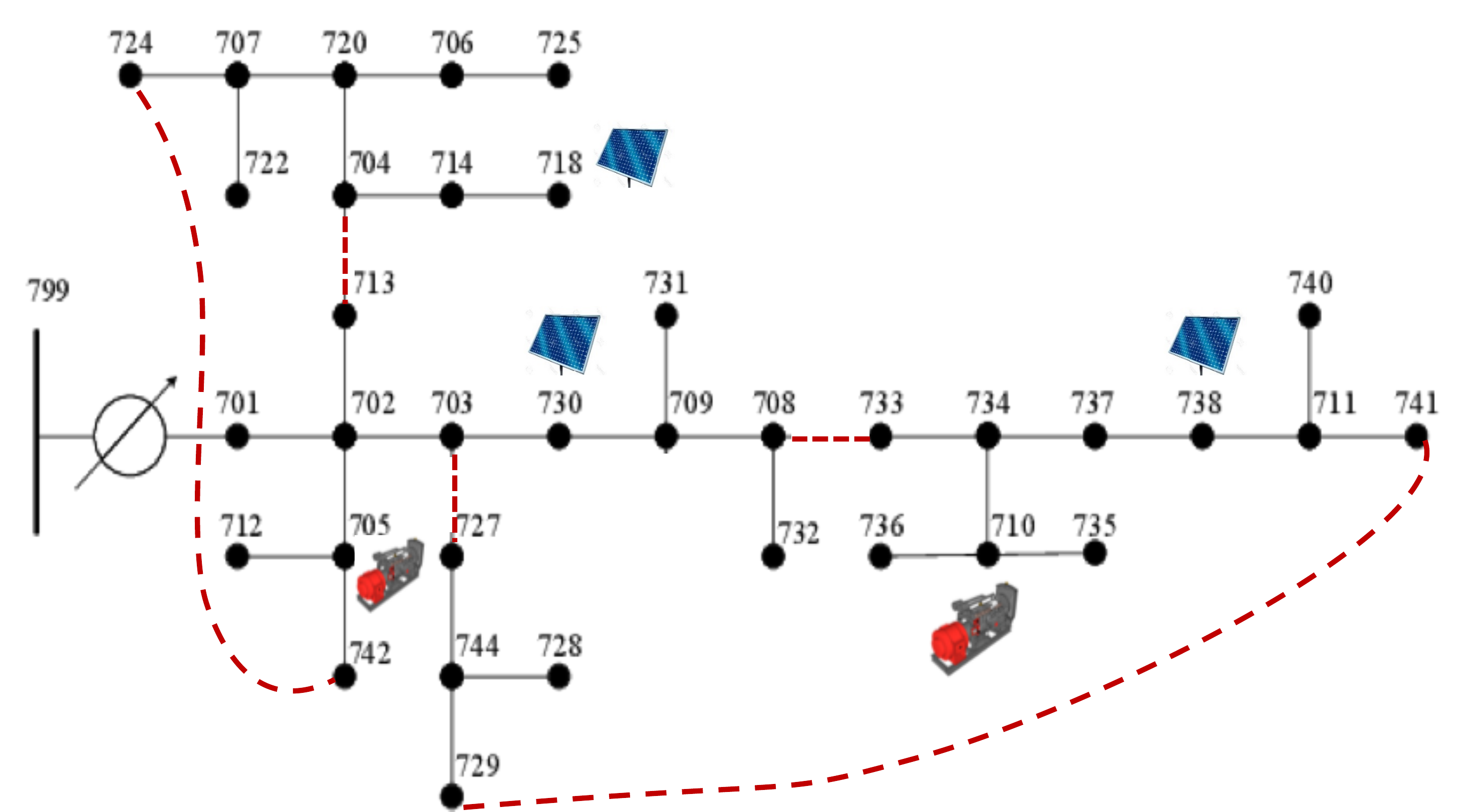}
	\caption{A modified IEEE 37-bus feeder showing existing lines and generators.}
	\label{fig:pds}
\end{figure}

To capture which buses are energized, introduce the vector of bus statuses $\bx\in\{0,1\}^{|\mcN|}$. Its $i$-th entry equals $1$ if bus $i$ is energized, and $0$ otherwise. If $v_i$ is the squared voltage magnitude and $p_i+jq_i$ the complex power injection on bus $i$, we enforce the constraints
\begin{subequations}\label{eq:limitsB}
	\begin{align}
	&x_i\in\{0,1\},\quad && \forall~ i\in\mcN\\
	&x_i\underline{v}_i\leq v_i\leq x_i\bar{v}_i,\quad &&\forall~i\in\mcN\label{seq:vol_lim}\\
	&x_i\underline{p}_i\leq p_i\leq x_i\bar{p}_i,\quad &&\forall~i\in\mcN\label{seq:p_lim}\\
	&x_i\underline{q}_i\leq q_i\leq x_i\bar{q}_i,\quad &&\forall~i\in\mcN.\label{seq:q_lim}
	\end{align}
\end{subequations}
Constraint \eqref{seq:vol_lim} ensures that voltages remain within voltage regulation limits for energized buses (e.g., $\pm3\%$ per unit); and set voltages to zero for non-energized buses. The substation voltage can be set by selecting $\underline{v}_0=\bar{v}_0=v_0$. Constraints \eqref{seq:p_lim}--\eqref{seq:q_lim} limit the complex power injections at energized nodes. The signed values for $\{\underline{p}_i,\overline{p}_i,\underline{q}_i,\overline{q}_i\}$ determine whether an injection corresponds to a generator; a fixed load with lagging or leading power factor; or an elastic load.

Not all DGs have black-start capabilities~\cite{lopes2007LVMG}: For instance, rooftop solar generators without energy storage can generate power only if the corresponding bus is already energized. On the other hand, diesel and gas-fired generators may feature black-start and grid-forming capabilities. To capture this functionality, define $\mcS_B\subseteq \mcN$ as the set of buses hosting black-start generators, and $\mcS_{NB}\subseteq \mcN$ as the set of buses with non-black-start generators. 

\subsection{Edge Variables and Constraints}\label{subsec:edges}
Let us now partition the set of edges $\mcE$ into:
\begin{itemize}
\item the subset $\mcE_O$ of out-of-service lines;
\item the subset $\mcE_I$ of in-service lines;
\item the subset $\mcE_S$ of switches; and 
\item the subset $\mcE_R$ of in-service voltage regulators.
\end{itemize}
Since non-remotely controlled switches cannot participate in DSR, they are handled as lines and belong to $\mcE_I$. Figure~\ref{fig:pds} depicts a feeder hosting $3$ solar generators; $2$ black-start diesel generators; $5$ switches; and $1$ voltage regulator.

Similar to buses, the vector of edge statuses $\by\in\{0,1\}^{|\mcE|}$ indicates which edges are closed. Vector $\by$ should satisfy
\begin{subequations}\label{eq:ybin}
	\begin{align}
&y_e\in\{0,1\},\quad && \forall~e\in\mcE_S\\
&y_e=0,\quad && \forall~ e\in\mcE_O\\
&y_e=1,\quad &&\forall~ e\in\mcE_I\cup\mcE_R\\
&y_e\underline{P}_e\leq P_e\leq y_e\bar{P}_e,\quad &&\forall~ e\in\mcE\label{seq:P_lim}\\
&y_e\underline{Q}_e\leq Q_e\leq y_e\bar{Q}_e,\quad &&\forall~ e\in\mcE\label{seq:Q_lim}
	\end{align}
\end{subequations}
 where $P_e+jQ_e$ is the complex flow on line $e$. Power flow limits are typically set as $\underline{P}_e=-\bar{P}_e$ and $\underline{Q}_e=-\bar{Q}_e$. 
 
Even though apparent flow limits of the form $P_e^2+Q_e^2\leq S_e^2$ can be added to our formulation, they result in a mixed-integer quadratic program, which does not scale as gracefully as an MILP. Alternatively, apparent constraints on line flows and bus injections can be handled by a polytopic inner or outer approximation of $P_e^2+Q_e^2\leq S_e^2$; see e.g., \cite{Jabr18}. This approach is not adopted here to keep the formulation uncluttered. 

\subsection{Voltage Drops and Regulators}\label{subsec:voltage}
To relate power injections, power flows, and voltages, we adopt the \emph{linearized distribution flow} (LDF) model~\cite{BW1}--\cite{BW2}. Albeit approximate, the LDF model has been engaged in various grid optimization tasks with satisfactory accuracy~\cite{VKZG16}. Given the complexity and uncertainty involved in DSR, the approximation error incurred by LDF becomes inconsequential. 

Upon ignoring ohmic losses on lines, the LDF model expresses bus injections as
\begin{subequations}\label{eq:lindist}
	\begin{align}
	&p_i = \sum_{e:(i,j) \in \mcE}P_e - \sum_{e:(j,i) \in \mcE}P_e,\quad \forall~i\in\mcN\label{seq:active}\\
	&q_i = \sum_{e:(i,j) \in \mcE}Q_e - \sum_{e:(j,i) \in \mcE}Q_e,\quad \forall~i\in\mcN.\label{seq:reactive}
		\end{align}
\end{subequations}
Constraint \eqref{eq:lindist} essentially imposes the complex power balance at each node. For node $i=0$, it implies that the injected power equals the total power withdrawn by the feeder. 

According to the LDF model, the voltage drop along line $e$ with impedance $r_e+jx_e$ can be approximated as
\begin{equation}\label{eq:vdrop}
y_e\left(v_i-v_j-2r_eP_e-2x_eQ_e\right)=0,~\forall~e:(i,j)\in\mcE\setminus\mcE_R.
\end{equation}
The drop occurs only if line $e$ is closed $(y_e=1)$; and \eqref{eq:vdrop} is not enforced for voltage regulators. 

Constraint \eqref{eq:vdrop} involves the bilinear terms $y_ev_i$, $y_ev_j$, $y_eP_e$, and $y_eQ_e$. As discussed in Section~\ref{sec:pre}, each one of these products can be replaced by an auxiliary variable, so that \eqref{eq:vdrop} can be posed as a linear equality constraint relating the four auxiliary variables. The auxiliary variable associated with $y_ev_i$ is related to $y_e$ and $v_i$ via four linear inequalities as in \eqref{eq:MC}. The same holds for the other three bilinear terms. Luckily, the McCormick linearization is used scarcely, since only a few edges represent switches.

We proceed with voltage regulators, which are modeled as ideal. This is without loss of generality since the impedance of a non-ideal regulator can be modeled as a line connected in series with the ideal regulator. A regulator can scale its secondary-side voltage by $\pm10\%$ by increments of $0.625\%$ using tap positions~\cite{hiskens2016tap}. The taps can be changed either remotely, or based on some automated control usually based on local voltage (and current) measurements. Due to space limitations, all regulators are assumed to be remotely controlled. 

Consider regulator $r:(i,j)\in\mcE_R$. Its voltage transformation ratio can be set to $1+0.00625\cdot t_r$, where $t_r\in\{0,\pm1,\ldots,\pm 16\}$ is its tap position. The transformation in terms of squared voltage magnitudes is
\begin{equation}\label{eq:VR}
v_j=(1+0.00625\cdot t_r)^2v_i.
\end{equation}
The quadratic dependence on $t_r$ is often replaced by a linear approximation~\cite{lopez2015tap}. Waiving this approximation, we pursue a simple yet \emph{exact} regulator model: The term in the parenthesis of \eqref{eq:VR} can take one out of $33$ possible values. These values are collected in vector $\bc\in\mathbb{R}^{33}$ whose $k$-th entry is $c_k:=\left[1+0.00625\cdot\left(k-17\right)\right]^2$. Vector $\bc$ is known beforehand and is common for all regulators. By introducing the tap status vectors $\bt_r$, the operation of regulators is modeled as
\begin{subequations}\label{eq:tap}
	\begin{align}
	&v_j=v_i \cdot \bt_r^\top\bc ,\quad &&\forall~r:(i,j)\in\mcE_R\label{seq:tap}\\
	&\bt_r\in\{0,1\}^{33},\quad &&\forall~r:(i,j)\in\mcE_R\label{seq:tap2}\\
	&\bt_r^\top\bone=1,\quad &&\forall~r:(i,j)\in\mcE_R.\label{seq:tap3}
	\end{align}
\end{subequations}
There is only one non-zero entry in $\bt_r$ due to \eqref{seq:tap3}. The bilinear product in \eqref{seq:tap} can be handled via the McCormick scheme. 


\subsection{Topological Constraints}\label{subsec:topology}
To deal with network constraints, let us introduce indicator vectors for paths and cycles. Because $\mcG$ is connected, there exists at least one path for each pair of nodes $i$ and $j$. For path $\mcP$, define its indicator vector $\bpi^{\mcP}\in\{0,1\}^{|\mcE|}$, such that $\pi^\mcP_e=1$ if $e\in\mcP$, and $\pi^\mcP_e=0$ otherwise. In essence, vector $\bpi^{\mcP}$ indicates which edges comprise $\mcP$, regardless their directionality. Similarly, for any cycle $\mcC$ in $\mcG$, define the cycle indicator vector $\bn^{\mcC}\in\{0,1\}^{|\mcE|}$, such that $n^\mcC_e=1$ if $e\in\mcC$, and $n^\mcC_e=0$ otherwise.

Distribution grids are typically operated in a tree (radial) structure to ease protection coordination. To enforce radiality, previous works were confined to unidirectional power flows~\cite{ccliu2014spanning}. Anticipating increasing penetration of renewables, we facilitate radiality even with reverse flows. To avoid the formation of cycles, we add the constraint
\begin{equation}\label{eq:tree}
\by^\top\bn^\mcC\leq \mathbf{1}^\top\bn^\mcC-1,\quad\forall\mcC
\end{equation}
for all cycles $\mcC$ in $\mcG$. Constraint \eqref{eq:tree} limits the number of closed edges along cycle $\mcC$ to be less than the total number of edges in $\mcC$. Due to the limited number of switches, there are few cycles $\mcC$. For instance, the system of Figure~\ref{fig:pds} has five switches giving rise to two cycles. 

In addition, if line $e:(i,j)\in\mcE$ is closed, the buses $i$ and $j$ must share the same status (both energized or not)
\begin{equation}\label{eq:xy}
|x_i-x_j|\leq 1-y_e,\quad\forall~e:(i,j)\in\mcE.
\end{equation}

\subsection{Coordinating Generators}\label{subsec:coordination}
Non-black-start generators (e.g., rooftop photovoltaics) can generate power only when they are connected to a substation or a running black-start generator; see~\cite{jianhui2018sequential}, \cite{lopes2007LVMG}. If there exists such path for generator $i\in\mcS_{NB}$, then constraint \eqref{eq:xy} implies $x_i=1$. Otherwise, the status $x_i=0$ must be enforced explicitly. To this end, identify all paths from generator $i\in \mcS_{NB}$ to all $j\in\mcS_B$. These paths are denoted by $\mcP_{i,k}$ for $k=1,\ldots,K_i$ and all $i\in \mcS_{NB}$; there are $21$ such paths in the feeder of Fig.~\ref{fig:pds}. For path $\mcP_{i,k}$, let $\bpi_{i,k}$ be its indicator vector and introduce the binary variable $\delta_{i,k}$ for which
\begin{subequations}\label{eq:nbs-delta}
	\begin{align}
&\delta_{i,k}\in\{0,1\},\quad && \forall\mcP_{i,k} \label{eq:nbs-delta:1}\\
&\delta_{i,k}\cdot \bone^\top\bpi_{i,k}\leq \by^\top\bpi_{i,k},\quad && \forall\mcP_{i,k}\label{eq:nbs-delta:2}\\
&x_i\leq\sum_{k}\delta_{i,k},\quad && \forall i\in\mcS_{NB}.\label{eq:nbs}
\end{align}
\end{subequations} 
By definition of the indicator vector $\bpi_{i,k}$, we have that $\by^\top\bpi_{i,k}\leq \bone^\top\bpi_{i,k}$ with equality only if path $\mcP_{i,k}$ is energized. If $\mcP_{i,k}$ is \emph{not} energized, then \eqref{eq:nbs-delta:1}--\eqref{eq:nbs-delta:2} imply $\delta_{i,k}=0$. If $\mcP_{i,k}$ is energized, then  $\delta_{i,k}$ can be either $0$ or $1$; yet bus $i$ is guaranteed to be energized by applying \eqref{eq:xy} along $\mcP_{i,k}$. Constraint \eqref{eq:nbs} entails that for bus $i$ to be energized, at least one of the paths $\{\mcP_{i,k}\}_k$ is energized. In other words, each non-black-start DG can run only if it is connected to a running black-start DG or a substation. 

When an island includes multiple black-start DGs, a simple coordination scheme is that the largest DG $i$ operates in PV mode $(v_i=v_0)$, and the rest in PQ mode~\cite{lopes2007LVMG}. Moreover, if black-start DGs operate in grid-connected mode, the substation should be treated as the largest generator and all DGs operate in PQ mode. To model this, identify all paths between each $i\in\mcS_B$ and all generators $j\in\mcS_B$ of larger rating. These paths are denoted by $\mcL_{i,\ell}$ and indexed by $\ell=1,\ldots,L_i$ for $i\in\mcS_B$; there are $8$ such paths in the feeder of Fig.~\ref{fig:pds}. The indicator vector for path $\mcL_{i,\ell}$ is denoted by $\bdell_{i,l}$. For each $\mcL_{i,\ell}$, introduce the  variable $\epsilon_{i,\ell}$, which equals $1$ if $\mcL_{i,\ell}$ is energized; and $0$, otherwise. With the help of $\epsilon_{i,\ell}$'s, the coordination of generators can be captured by
\begin{subequations}\label{eq:bs-delta}
	\begin{align}
	&\epsilon_{i,\ell}\in\{0,1\},\quad &&\forall \mcL_{i,\ell}\\
	& (\by-\bone)^\top\bdell_{i,l}+1\leq\epsilon_{i,\ell} \leq \frac{\by^\top\bdell_{i,l}}{\bone^\top\bdell_{i,l}},\quad &&\forall \mcL_{i,\ell}.
	\end{align}
\end{subequations}
If path $\mcL_{i,\ell}$ is energized, then $\by^\top\bdell_{i,l}=\bone^\top\bdell_{i,l}$, and \eqref{eq:bs-delta} entails $\epsilon_{i,\ell}=1$. If $\mcL_{i,\ell}$ is not energized, then $\by^\top\bdell_{i,l}<\bone^\top\bdell_{i,l}$. Because $\by^\top\bdell_{i,l}$ counts the number of closed lines in $\mcL_{i,\ell}$, it holds $\by^\top\bdell_{i,l}\leq\bone^\top\bdell_{i,l}-1$ and so \eqref{eq:bs-delta} entails $\epsilon_{i,\ell}=0$.

To model the operation mode for generator $i\in\mcS_B$, introduce variable $\epsilon_i$ that equals $1$ if the generator operates in PQ, and $0$ when in PV mode. Using $\epsilon_{i,\ell}'s$, the coordination of generator modes is accomplished as
\begin{subequations}\label{eq:bs}
	\begin{align}
	&\epsilon_i \in\{0,1\},\quad && \forall i\in\mcS_B\label{seq:bs-b}\\
		&\max_{\ell}\epsilon_{i,\ell}\leq \epsilon_i\leq\sum_{\ell}\epsilon_{i,
		\ell},\quad && \forall i\in\mcS_B\label{seq:bs-a}\\
		&|v_i-v_0|\leq \epsilon_i v_0,\quad && \forall i\in\mcS_B.\label{seq:bs-c}
	\end{align}
\end{subequations}
Constraint \eqref{seq:bs-a} ensures that if any of the paths $\{\mcL_{i,\ell}\}_{\ell}$ is energized, then $\epsilon_i=1$ and \eqref{seq:bs-c} follows trivially due to voltage regulation. This scenario means that generator $i$ is connected to a larger black-start generator or the substation, and hence operates in PQ mode. On the other hand, if $\epsilon_{i,\ell}=0$ for all $\ell$, generator $i$ is the largest on the island. In this case, constraint \eqref{seq:bs-a} yields $\epsilon_i=0$, and \eqref{seq:bs-c} sets $v_i=v_0$.

\subsection{Objective Function}\label{subsec:objective}
Let vectors $\bx^0$ and $\by^0$ represent the post-outage statuses of nodes and lines. A meaningful restoration objective is to find a grid topology $\by$ and generation dispatch that maximize the total served load. Among several restoration schemes, an operator may prefer the schemes with fewer line switching operations. Moreover, a usual practice dictates that the restoration process must not de-energize an already energized bus. The DSR problem can be now formulated as
\begin{align}\label{eq:DSR1}
\min~&~\sum_{i\in\mcN\setminus(\mcS_B\cup\mcS_{NB}\cup\{0\})}p_i+\lambda\mathbf{1}^\top|\by-\by^0|\\
\mathrm{s.to}~&~ \eqref{eq:xy}-\eqref{eq:bs}, \bx\geq\bx^0
\end{align}
where parameter $\lambda\geq0$ quantifies the importance of restoration schemes with fewer switching operations. Setting $\lambda=0$ yields the scheme with the maximum load served. Problem \eqref{eq:DSR1} is an MILP and can be solved by off-the-shelf solvers.


\section{Numerical Tests}\label{sec:tests}
The developed DSR approach was tested on a modified version of the IEEE 37-node feeder converted to its single-phase equivalent~\cite{guido2018analytics}; see Figure~\ref{fig:pds}. Two black-start DGs of capacities $459.3$ and $918.5$~kW were placed on nodes $705$ and $710$, respectively. Three non-black-start DGs were placed on buses $718$, $730$, and $738$ with capacities set to half the load on the associated buses. The switchable lines include $3$ existing and two additional lines shown as dashed. The (re)active loads on buses $701$, $722$, $737$, and $738$ were elastic with their minimum set to half the nominal bus load. All tests were run using MATLAB-based toolbox YALMIP along with the mixed-integer solver CPLEX~\cite{YALMIP},~\cite{cplex}; on a $2.7$ GHz Intel Core i5 computer with 8 GB RAM; and for $\lambda=10^{-3}$.
 
\begin{figure}[t]
	\centering
	\includegraphics[scale=0.28]{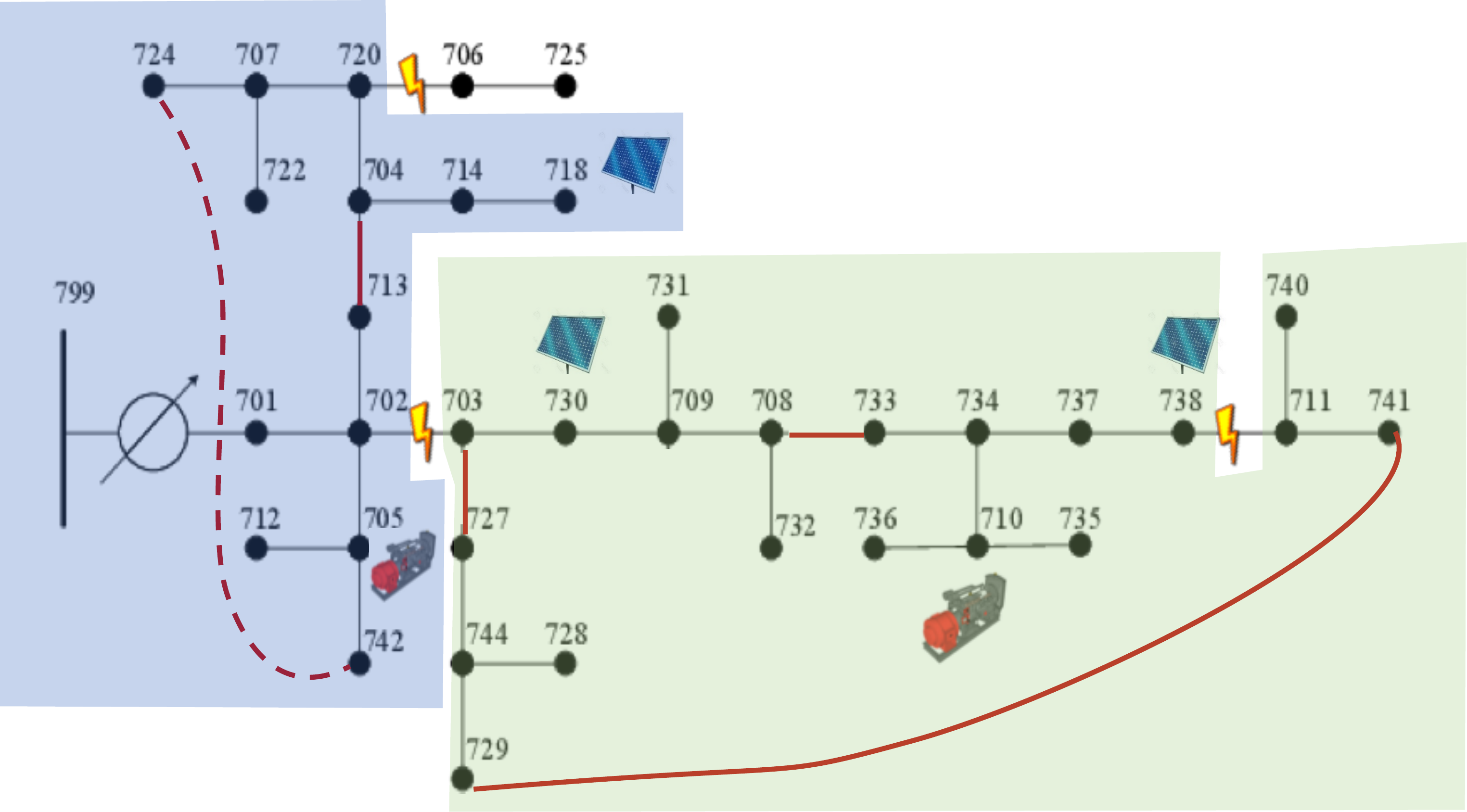}
	\caption{The feeder of Figure~\ref{fig:pds} restored after a $3$-line outage.}
	\label{fig:pds-re}
\end{figure}

\begin{figure}[t]
	\centering
	\includegraphics[scale=0.22]{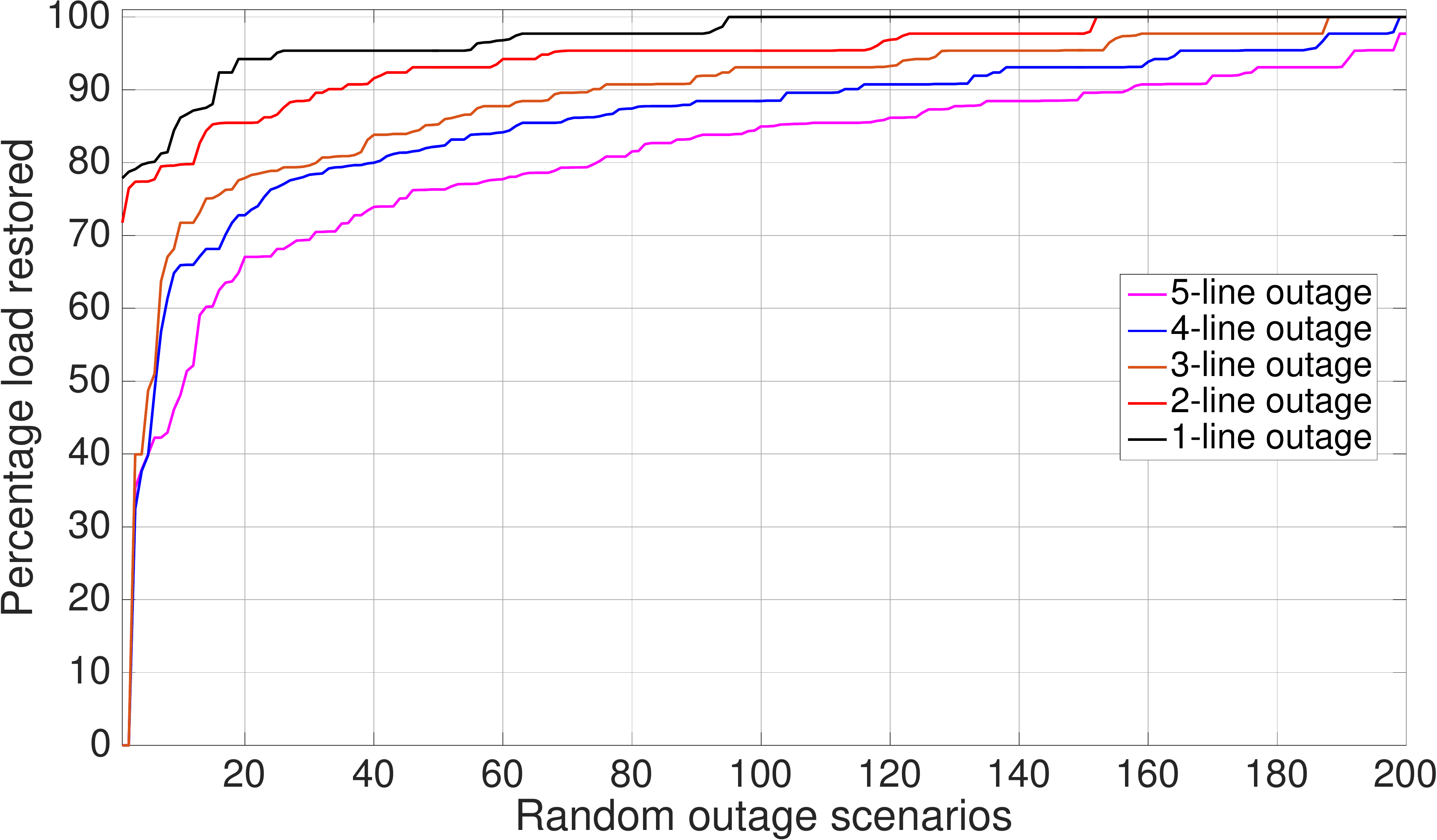}
	\caption{The (ordered) percentage of load restored after $1$--$5$ line outages.}
	\label{fig:load_res}
\end{figure}

The performance of our DSR scheme was tested for an outage scenario of three line outages shown in Figure~\ref{fig:pds-re}. The restored system comprises of two sub-networks, while buses $706$ and $725$ could not be restored. The DG on bus $710$ serves as the reference bus for the second island. 

The computational performance of the MILP in \eqref{eq:DSR1} was tested using $1,000$ random outage scenarios, $200$ scenarios for each number of $1-5$ lines in outage. The maximum available solar generations were drawn from a uniform distribution based on the respective rated sizes. The running times for solving \eqref{eq:DSR1} reported in Table~\ref{tbl:time} demonstrate that our DSR scales well for single- and multiple-line outages alike. The percentage load restored for the various line outages is shown (ordered) in Figure~\ref{fig:load_res}. As anticipated, the total load restored decreases as the number of outages increases.

\begin{table}
\renewcommand{\arraystretch}{1}
\caption{Running Times for the MILP of \eqref{eq:DSR1}}
\vspace*{-1em}
\label{tbl:time} \centering
\begin{tabular}{|l|r|r|r|r|r|}
\hline
\hline
Number of outaged lines & $1$ & $2$ & $3$ & $4$ & $5$ \\
\hline
\hline
Maximum running time [sec] & 1.04 & 0.96 & 2.69 & 3.96 & 1.77\\
Median running time [sec] & 0.79 & 0.78 & 0.77 & 0.81 & 0.73\\
\hline
\hline
\end{tabular}
\end{table}

\section{Conclusions}\label{sec:conclusions}
The developed DSR scheme features optimal formation of islands; incorporates voltage regulators; allows for multiple DGs on each island and establishes a coordination hierarchy amongst them. Numerical tests demonstrate the correctness of the MILP formulation and that its complexity scales well in moderately-sized feeders. Its scalability can be attributed to three key points: \emph{i)} the unique use of indicator vectors for cycles and paths over the infrastructure graph; \emph{ii)} the McCormick linearization; and \emph{iii)} the approximate LDF model. Although framed within the DSR paradigm, this work sets the solid foundations for several grid optimization tasks including reconfiguration for power loss minimization and Volt/VAR control. We are currently working towards extending this scheme to its multi-step variant; incorporating switched capacitor banks and locally-controlled voltage regulators; and considering unbalanced multi-phase feeders. 


\balance
\bibliography{myabrv,power}
\bibliographystyle{IEEEtran}
\end{document}